\documentclass[11pt,reqno,tbtags]{article}
\usepackage{amsmath}
\usepackage{amsthm}

\numberwithin{equation}{section}

\newtheorem{theorem}{Theorem}[section]
\newtheorem{corollary}[theorem]{Corollary}
\newtheorem{lemma}[theorem]{Lemma}
\newtheorem{proposition}[theorem]{Proposition}
\newtheorem{claim}[theorem]{Claim}
\newtheorem{example}[theorem]{\sl Example}
\newtheorem{algorithm}[theorem]{Algorithm}

\theoremstyle{definition}
\newtheorem{remark}[theorem]{Remark}

\newcommand{\lc}{\left\lceil}
\newcommand{\lf}{\left\lfloor}
\newcommand{\rc}{\right\rceil}
\newcommand{\rf}{\right\rfloor}

\newcommand{\CC}{{\bf  C}}

\newcommand{\PP}{{\bf  P}}
\newcommand{\QQ}{{\bf  Q}}

\newcommand{\XX}{{\bf  X}}

\newcommand{\pp}{{\bf  p}}

\newcommand{\xx}{{\bf  x}}
\newcommand{\yy}{{\bf  y}}
\newcommand{\zz}{{\bf  z}}

\newcommand{\mt}{{\tilde{m}}}

\newcommand{\Tc}{{\cal T}}

\newcommand{\Xc}{{\cal X}}
\newcommand{\Zc}{{\cal Z}}

\newcommand{\begp}{\begin{proposition}}
\newcommand{\enp}{\end{proposition}}
\newcommand{\begt}{\begin{theorem}}
\newcommand{\ent}{\end{theorem}}
\newcommand{\begl}{\begin{lemma}}
\newcommand{\enl}{\end{lemma}}
\newcommand{\begc}{\begin{corollary}}
\newcommand{\enc}{\end{corollary}}
\newcommand{\begcl}{\begin{claim}}
\newcommand{\encl}{\end{claim}}
\newcommand{\begr}{\begin{remark}}
\newcommand{\enr}{\end{remark}}
\newcommand{\begal}{\begin{algorithm}}
\newcommand{\enal}{\end{algorithm}}
\newcommand{\begd}{\begin{definition}}
\newcommand{\enf}{\end{definition}}
\newcommand{\begx}{\begin{example}}
\newcommand{\enx}{\end{example}}
\newcommand{\bega}{\begin{array}}
\newcommand{\ena}{\end{array}}
\newcommand{\s}{\hspace*{5ex}}

\newcommand{\Repeat}{{\tt repeat\ }}
\newcommand{\Until}{{\tt until\ }}
\newcommand{\For}{{\tt for\ }}
\newcommand{\To}{{\tt \ to\ }}

\newcommand{\Return}{{\tt return\ }}

\newcommand{\sset}{\leftarrow}

\newcommand{\Random}{\mbox{\tt Random()}}
\newenvironment{code}{\begin{list}{\hspace*{0em}}{}}
{\end{list}}

\def\rompar(#1){\textup(#1\textup)}    % usage: \rompar(...)

\newcommand\gk{\kappa}

\newcommand\PPbar{{\overline \PP}}

\newcommand{\refS}[1]{Section~\ref{#1}}
\newcommand{\refT}[1]{Theorem~\ref{#1}}

\newcommand{\refL}[1]{Lemma~\ref{#1}}
\newcommand{\refP}[1]{Proposition~\ref{#1}}
\newcommand{\refA}[1]{Algorithm~\ref{#1}}
\newcommand{\refR}[1]{Remark~\ref{#1}}

\newcommand\ie{i.e.\spacefactor=1000}
\newcommand\eg{e.g.\spacefactor=1000}

\newcommand\cf{{cf.}\spacefactor=1000}
\newcommand\as{a.s.\spacefactor=1000}

\newcommand\nopf{\qed}   % for theorem without proof
\begin{document}

\setcounter{page}{0}
\thispagestyle{empty}

\begin{center}
{\Large \bf Interruptible Exact Sampling in the Passive Case \\}
%\vspace{.15in}
\normalsize

\vspace{4ex}
{\sc Keith Crank \footnotemark} \\
\vspace{.1in}
Division of Mathematical Sciences \\
\vspace{.1in}
National Science Foundation \\
\vspace{.1in}
{\tt kcrank@nsf.gov} \\
\vspace{.2in}
{\sc and} \\
\vspace{.1in}
{\sc James Allen Fill \footnotemark} \\
\vspace{.1in}
Department of Mathematical Sciences \\
\vspace{.1in}
The Johns Hopkins University \\
\vspace{.1in}
{\tt jimfill@jhu.edu} and {\tt http://www.mts.jhu.edu/\~{}fill/} \\
\end{center}
\vspace{3ex}

\begin{center}
{\sl ABSTRACT} \\
\end{center}

We establish, for various scenarios, whether or not interruptible exact stationary
sampling is possible when a finite-state Markov chain can only be viewed passively. 
In particular, we prove that such sampling is not possible using a single copy of
the chain.  Such sampling is possible when enough copies of the chain are available,
and we provide an algorithm that terminates with probability one.
%\vspace{.05in}
\bigskip
\bigskip

\begin{small}

\par\noindent
\emph{AMS} 2000 \emph{subject classifications.}  Primary 60J10, 68U20;
secondary 60G40, 62D05, 12D05.
%62E25
%not yet: 60-04
%
%Here
%
%60xxx = Probability theory and stochastic processes
%not yet: %60-04 = Explicit machine computation and programs
%60Jxx = Markov processes
%60J10 = Markov chains with discrete parameter
%60Gxx = Stochastic processes
%60G40 = Stopping times; optimal stopping times; gambling theory
%68xxx = Computer science
%68Uxx = Computing methodologies and applications
%68U20 = Simulation
%68Pxx = Theory of data
%68P05 = Data structures
%68P10 = Searching and sorting
%62xxx = Statistics
%62D05 = Sampling theory; sample surveys
%62E25 = Monte Carlo studies
%65-XX = Numerical analysis
%65Cxx = Probabilistic methods, simulation and stochastic differential equations
%65C05 = Monte Carlo methods
%65C10 = Random number generation
%65C40 = Computational Markov chains
%12-XX = Field theory and polynomials 
%12Dxx = Real and complex fields
%12D05 = Polynomials: factorization
%05C05 = Trees
\medskip
\par\noindent
\emph{Date.} February~14, 2002.
\end{small}

\footnotetext[1]{The views and findings presented here are those of the author and
do not necessarily reflect the opinion of the National Science Foundation.}

\footnotetext[2]{This author's research was supported by NSF grants DMS--9803780 and
DMS--0104167, and by The Johns Hopkins University's Acheson J.~Duncan Fund for the
Advancement of Research in Statistics; and was carried out in part while
author was Visiting Researcher, Theory Group, Microsoft Research.}

\newpage
\addtolength{\topmargin}{+0.5in}
%%%THIS IS THE END OF STANDARD PREAMBLE-TYPE STUFF%%%
%%%%%%%%%%%%%%%%%%%%%%%%%%%%%%%%%%%%%%%%%%%%%%%%%%%%%%%
%%%%%%%%%%%%%%%%%%%%%%%%%%%%%%%%%%%%%%%%%%%%%%%%%%%%%%%

\section{Introduction and summary}
\label{S:intro}

In recent years a large number of articles have been written about exact sampling
(also called perfect sampling) using Markov chains.  See~\cite{Wilson} for an
overview.  The rough idea is as follows.  One wishes to sample from the unique
stationary distribution~$\pi$ of an observed irreducible Markov chain.  At each
transition of the chain, a decision is made whether to continue observing the chain
or to stop.  When the observation is stopped, a value~$S$ is output and it is
desired that, for all states~$i$, $P(S = i\,|\,\mbox{one ever stops observing the
chain}) = \pi_i$.  The decision about whether to stop at a particular time is made
on the basis of the evolution of the chain up through that time, possibly together
with some additional randomness independent of the chain.

The goal of our research leading to this paper
was to determine whether or not it is possible to carry out interruptible exact
sampling for finite-state chains in what Propp and Wilson~\cite{PW98} call the
passive setting.  (We will explain in \refS{S:passive} what is meant by
``interruptible'' and ``the passive setting''.)  Our central result is the
following:
\begin{center}
\emph{Interruptible exact sampling is not possible when one observes only a single
trajectory.}
\end{center}
This result remains true even if we assume that the chain is aperiodic and
reversible.  [See \refR{R:rev}(b).] However, interruptible exact sampling \emph{is}
possible for an $N$-state chain when one is able to observe, simultaneously, $N$
trajectories.  Here is a guide to our specific results.
\begin{enumerate}

\item (positive:)\ \ We provide an algorithm (\refA{A:algogen}) which, given an
irreducible Markov chain on~$N$ states as input, produces in (random) finite time
an exact sample from the tree distribution, and hence also an exact sample
from~$\pi$.  (The tree distribution is defined in \refS{S:tree}.)   The algorithm
is interruptible, but requires~$N$ independent synchronized trajectories from the
chain.  (See \refT{T:algogen}.)

\item (negative:)\ \ There
is no algorithm in the passive setting for obtaining an
observation from the stationary distribution of an irreducible aperiodic Markov
chain on~$N$ states which uses fewer than~$N$ independent trajectories from the
chain and which is both interruptible and exact.  (See \refT{T:neg1}.)

\item (negative:)\ \ There is no algorithm in the passive setting for obtaining an
observation from the common stationary distribution of any finite number of
independent irreducible aperiodic Markov chains on~$N$ states (with possibly
different transition matrices) which is both interruptible and exact.  (See
\refT{T:neg2}.)  This remains true even if we assume that all of the chains are
reversible.  [See \refR{R:rev}(a).]

\end{enumerate}

\section{Background}
\label{S:back}

In 1992, Asmussen, Glynn, and Thorisson~\cite{AGT} demonstrated that exact sampling
from a Markov chain is possible under certain circumstances.  They also proved that
it is not possible to obtain an exact sample from an arbitrary Markov chain without
some prior knowledge about the chain; in particular, the size of the state space
must be known.  Although their paper does provide a method for generating exact
samples from an $N$-state Markov chain when~$N$ is known, the paper is primarily of
a theoretical nature, and the method is complicated and inefficient.

In 1995, Lov\'asz and Winkler~\cite{LW} provided a simpler and more efficient
algorithm for obtaining an exact sample from an irreducible $N$-state Markov
chain.  Although not mentioned explicitly in their paper, the method described in
Section~3 of Lov\'asz and Winkler can in fact be used to obtain an exact sample
from the tree distribution of the Markov chain (as defined in \refS{S:tree}).
Aldous~\cite{Aldous}, Broder~\cite{Broder}, and Propp and Wilson~\cite{PW98} also
describe algorithms for sampling from the tree distribution.  Propp and
Wilson~\cite{PW98} discuss and compare these and other methods of sampling from the
tree distribution, and from the stationary distribution.  Their discussion includes
consideration of such issues as whether or not the sampling is exact or
interruptible.  To our knowledge, the question of whether interruptible exact
sampling is possible in the passive case (as described in \refS{S:passive}) has not
previously been considered.

\section{Preliminaries}
\label{S:prel}

\subsection{The tree distribution}
\label{S:tree}

Throughout this paper we consider only finite-state irreducible Markov chains.  We
assume that the number of states, call it~$N$, is known; in fact, it turns out that
we may as well assume (and so we do) that the state space is known to be $[N] :=
\{1, \ldots, N\}$.  We denote the transition matrix of such a chain generically by
$\PP = (p_{ij})$.
%_{i \in [N], j \in [N]}$.

An irreducible Markov chain on $[N]$ can be viewed equivalently as a random walk on
a connected weighted directed graph~$G$.  The vertex set of~$G$ is $[N]$, and there
is an edge from~$i$ to~$j$, with weight $p_{ij}$, if and only if $p_{ij} > 0$. 

For the moment, let us consider an \emph{undirected} graph~$G$ with vertex set
$[N]$.  Then a subgraph~$T$ of~$G$ is called a \emph{spanning tree} if it contains
all~$N$ vertices and is connected and acyclic.  From any spanning tree, we obtain a
\emph{directed spanning tree} by assigning a direction to each edge.  A directed
spanning tree is called an \emph{arborescence} rooted at a given vertex~$r$ if all
edges are directed towards~$r$. 

We define the \emph{weight} $w(T)$ of an arborescence~$T$ with edges $\{e_l\}$ as
$w(T) := \prod_{l = 1}^{N - 1} p(e_l)$, where $p(e_l) := p_{ij}$ if $e_l$ is
directed from~$i$ to~$j$.  For the remainder of this paper, when we
say ``tree'' we mean an arborescence~$T$ with $w(T) > 0$.
The \emph{tree distribution} of the Markov chain is the probability distribution on
trees obtained by normalizing the weights $w(T)$ so as to sum to unity.

The Markov chain tree theorem is the well-known result (see, for example, \cite{LW}
or~\cite{AF}) that the stationary distribution~$\pi$ of the chain can be expressed
simply in terms of the tree distribution:
$$
\pi_i = w_i / w, \quad i \in [N],
$$
where, writing $\Tc_i$ for the set of trees rooted at~$i$ and~$\Tc$ for
$\cup_{i \in [N]} \Tc_i$,
$$
w_i := \sum_{T \in \Tc_i} w(T), \qquad w := \sum_{i \in [N]} w_i = \sum_{T \in
\Tc} w(T).
$$
In particular, any algorithm for sampling from the tree distribution provides a
means of sampling from~$\pi$:\ simply output the root of the tree.  

\subsection{The passive case; interruptible exact sampling}
\label{S:passive}

Propp and Wilson~\cite{PW98} distinguish between the active setting and the passive
setting for sampling using a Markov chain. In the active setting, an algorithm is
assumed to have access at all times to a transition generator, that is, to a routine
which, given any input state~$i$, generates an observation~$j$ from the probability
distribution $(p_{ij}: j \in [N])$, independent of all previously generated
observations.  In particular, a user can generate a trajectory from~$\PP$ with any
desired initial state.  In the passive setting, the algorithm has no control over
the initial state and can only watch passively as the chain transitions from one
state to the next.

We now explain what is meant by an (on-line, Markov-chain-based) interruptible exact
sampling algorithm in the passive case; for simplicity, we will do this explicitly
only in the case that a single trajectory from the chain is available and the
desired output is an observation from the stationary distribution~$\pi$ (rather
than one from the tree distribution).  Informally, an exact sampling algorithm must
take as input a trajectory from the given Markov chain; possibly using external
randomization to make its decisions, it watches the chain only until some finite
time and then returns an observation distributed according to~$\pi$.  
(Important note:\ The state returned is not necessarily the state of the chain at
the stopping time.)  More formally, we can define an exact sampling algorithm as a
collection of functions $\phi_{k,i}: [N]^{k + 1} \to [0, 1]$ [with $\phi_{k,i}(x_0,
\ldots, x_k)$ to be interpreted informally as the conditional probability that the
algorithm stops by time~$k$ and outputs~$i$, given that it sees the trajectory
$(x_0, \ldots, x_k)$ through time~$k$] having the following properties, where~(iii)
and~(iv) must hold for all~$\pi$, for all~$\rho$, and for all irreducible
transition matrices $\PP = (p_{ij})$ on $[N]$ with stationary distribution~$\pi$:  
\begin{enumerate}

\item $\forall k \geq 0\ \ \forall (x_0, \ldots, x_k) \in
[N]^{k + 1}:\ \ \sum_j \phi_{k,j}(x_0, \ldots, x_k) \leq 1$;

\item $\forall i \in [N]\ \ \forall k \geq 0\ \ \forall (x_0, x_1, \ldots) \in
[N]^{\infty}:\ \ \phi_{k,i}(x_0, \ldots, x_k) \uparrow\mbox{\ as\ }k \uparrow$;

\item $\lim_{k \uparrow \infty} \sum_{j \in [N]} \sum_{x_0, x_1, \ldots, x_k}
\rho_{x_0} p_{x_0, x_1} \cdots p_{x_{k - 1}, x_k}\,\phi_{k,j}(x_0, x_1, \ldots, x_k)
> 0$;

\item $\forall i \in [N]:\ \ \lim_{k \uparrow \infty} \frac{\sum_{x_0, x_1, \ldots,
x_k} \rho_{x_0} p_{x_0, x_1} \cdots p_{x_{k - 1}, x_k}\,\phi_{k,i}(x_0, x_1, \ldots,
x_k)}{\sum_{j \in [N]} \sum_{x_0, x_1, \ldots, x_k} \rho_{x_0} p_{x_0, x_1} \cdots
p_{x_{k - 1}, x_k}\,\phi_{k,j}(x_0, x_1, \ldots, x_k)} = \pi_i$. 

\end{enumerate}
In terms of the chain~$X$ observed and the stopping time~$\tau$ and output state~$S$
for the algorithm, the properties can be interpreted informally as (i) $P(\tau \leq
k\,|\,X) \leq 1$; (ii)~$P(\tau \leq k\,|\,X) \uparrow$ as $k \uparrow$; (iii)
$P(\tau < \infty) > 0$; and (iv)~$P(S = i\,|\,\tau < \infty) \equiv \pi_i$.  When
the strengthening
\begin{enumerate}

\item[(iii$'$)] $\lim_{k \uparrow \infty} \sum_{j \in [N]} \sum_{x_0, x_1, \ldots,
x_k} \rho_{x_0} p_{x_0, x_1} \cdots p_{x_{k - 1}, x_k}\,\phi_{k,j}(x_0, x_1, \ldots,
x_k) = 1$

\end{enumerate}
[interpreted as $P(\tau < \infty) = 1$] of~(iii) holds, we will call the algorithm
\emph{terminating}.  When~(iv) can be strengthened to
\begin{enumerate}

\item[(iv$'$)] $\forall i \in [N]\ \ \forall k \geq 0$:
\begin{eqnarray*}
\lefteqn{\sum_{x_0, x_1, \ldots, x_k} \rho_{x_0} p_{x_0, x_1} \cdots p_{x_{k - 1},
           x_k}\,\phi_{k,i}(x_0, x_1, \ldots, x_k)} \\
 &=& \pi_i \times \sum_{j \in [N]} \sum_{x_0, x_1, \ldots, x_k} \rho_{x_0} p_{x_0,
       x_1} \cdots p_{x_{k - 1}, x_k}\,\phi_{k,j}(x_0, x_1, \ldots, x_k)
\end{eqnarray*}

\end{enumerate}
[interpreted as the independence $P(\tau \leq k, S = i) \equiv P(\tau \leq k)
\pi_i$ of~$\tau$ and $S \sim \pi$], we say that the algorithm is
\emph{interruptible}.  An interruptible algorithm can be aborted without biasing
output; see the discussion in~\cite{Fill}.  For active-case algorithms, the
leading example of a non-interruptible algorithm is coupling from the
past~\cite{PW96}, while interruptible algorithms include cycle  popping~\cite{PW98},
Fill's rejection-based algorithm~\cite{Fill}~\cite{FMMR}, and the Randomness
Recycler~\cite{FH}.  The results of this paper, both positive and negative, are for
interruptible algorithms.

\section{A terminating algorithm for interruptible exact sampling in the passive
case}
\label{S:algo}

In this section we present a terminating algorithm for interruptible exact 
stationary sampling in the passive case, assuming that one can watch~$N$
synchronized copies $X_i = (X_i(t): t = 0, 1, \ldots)$, $i \in [N]$, of a Markov
chain with state space $[N]$ and irreducible transition matrix~$\PP$.  We allow
arbitrary initial distribution~$\rho$ for the $N$-variate chain $\XX := (X_1,
\ldots, X_N)$, but we assume that $X_1, \ldots, X_N$ are conditionally independent
given the initial state $(X_1(0), \ldots, X_N(0))$.  The algorithm will
produce an observation from the tree distribution corresponding to~$\PP$ (recall
\refS{S:tree}).   

\subsection{The algorithm in a restricted setting}
\label{S:restricted}

In
this subsection we present a terminating algorithm for interruptible exact
tree-sampling in the passive case that works under the following additional
restriction on~$\PP$:
$$
\mbox{{\bf Assumption~A:\ \ }$p_{1j}>0$ for all $j \in [N]$}.
$$
While this assumption may seem unreasonably restrictive, we will show in
\refS{S:general} how a simple modification of the algorithm can handle the more
general case.

To describe the algorithm we first define the following events for even $t \geq 2$:
\begin{eqnarray*}
A_t &:=& \cap_{l \in [N]} \{X_l(t - 2) = 1\}, \\
B_t &:=& \{X_1(t - 1) = 1\}, \\
C_t &:=& \{\{X_1(t), X_2(t - 1), \ldots, X_N(t - 1)\} = [N]\}, \\
D_t(T)
    &:=& \{\mbox{the graph with directed edges from $X_l(t - 1)$ to $X_l(t)$, $2
           \leq l \leq N$,} \\
    &  & \qquad \mbox{is the arborescence~$T$}\}, \\
D_t &:=& \cup_{T \in \Tc} D_t(T), \\
E_t(T)
    &:=& A_t \cap B_t \cap C_t \cap D_t(T), \\
E_t &:=& \cup_{T \in \Tc} E_t(T)
      =  A_t \cap B_t \cap C_t \cap D_t.
\end{eqnarray*}
\newpage  

\begal[Terminating interruptible tree-sampling, under Assumption~A]
\label{A:algoA}
{\rm For even $t \geq 2$, let $E_t(T)$ and $E_t$ be defined as above, and let $E_0
:= \emptyset$.  The algorithm is:
\begin{code}
\item \s   $t \sset 0$
\item \s   \Repeat
\item \s\s   $t \sset t + 2$
\item \s   \Until $E_t$ holds
\item \s   $S \sset T$, for the unique $T \in \Tc$ such that $E_t(T)$ holds
\item \s   \Return $S$
\end{code}
}
\enal

\begin{theorem}
\label{T:algoA}
When Asumption~A holds, \refA{A:algoA} is a terminating algorithm for interruptible
exact tree-sampling.
\end{theorem}

\begin{proof}
Let $\tau$ denote the supremum of the values of the variable~$t$ during the
operation of \refA{A:algoA}.  Now fix a candidate value~$t$ of~$\tau$.  Let~$T \in
\Tc$ be an arborescence, say with edges $e_l$ directed from~$i_l$ to~$j_l$, which
we choose to index (in some arbitrary but fixed order) by $l \in \{2, \ldots,
N\}$.  The event $D_t(T)$ is a disjoint union of $(N - 1)!$ subevents, with each
subevent corresponding to a way of mapping the $N - 1$ transitions $(X_l(t - 1),
X_l(t))$ to the $N - 1$ edges $e_l$.  These subevents will all enter symmetrically
into the calculation below of $P(\tau = t,\,S = T)$.  One such subevent is
$$
D'_t(T) := \cap_{l = 2}^N \{(X_l(t - 1), X_l(t)) = (i_l, j_l)\}.
$$
Let $\{i_1\}$ denote the singleton $[N] \setminus \{i_2, \ldots, i_N\}$.

Define
$$
A'_t := \{\tau \geq t - 2\} \cap A_t, \quad \mbox{even $t \geq 2$}.
$$
Then, using the Markov property and independence of the trajectories,
\begin{eqnarray*}
\lefteqn{P(\tau = t,\,S = T)} \\
 &=& P(\{\tau \geq t - 2\} \cap E_t(T)) 
  =  P(A'_t \cap B_t \cap C_t \cap D_t(T)) \\
 &=& (N - 1)!\,P(A'_t \cap B_t \cap C_t \cap D'_t(T)) \\
 &=& (N - 1)!\,P(A'_t)\,P(B_t\,|\,X_1(t - 2) = 1) \left[ \prod_{l = 2}^N P(X_l(t -
       1) = i_l\,|\,X_l(t - 2) = 1) \right] \\
 & & \qquad \times  P(X_1(t) = i_1\,|\,B_t) \left( \prod_{l = 2}^N p_{i_l, j_l}
       \right) \\
 &=& (N - 1)!\,P(A'_t)\,p_{11} \left[ \prod_{l = 2}^N p_{1, i_l} \right] p_{1, i_1}
       w(T)
  =  (N - 1)!\,P(A'_t)\,p_{11} \left( \prod_{l = 1}^N p_{1l} \right) w(T).
\end{eqnarray*}
Summing over $T \in \Tc$ we find
$$
P(\tau = t) = (N - 1)!\,P(A'_t)\,p_{11} \left( \prod_{l = 1}^N p_{1l} \right) w
$$
and therefore
$$
P(\tau = t,\,S = T) = P(\tau = t) \frac{w(T)}{w},
$$
which shows that \refA{A:algoA} is an interruptible exact tree-sampling algorithm. 
Using the fact that~$\XX$ visits $(1, \ldots, 1)$ at even times infinitely often
(\as) together with the strong Markov property of~$\XX$, it is clear that
termination occurs at the first success in an almost surely infinite sequence of
Bernoulli trials with success probability $p_{11} \left( \prod_{l = 1}^N p_{1l}
\right) w > 0$ (note that this is where Assumption~A is used).  Thus $P(\tau <
\infty) = 1$, that is, \refA{A:algoA} is terminating.
\end{proof}

\subsection{The algorithm in the general setting}
\label{S:general}

To avoid needing Assumption~A, we can use the averaging technique of Lov\'asz and
Winkler~\cite{LW}.  Let $\PP^k$ be the $k$-step transition matrix of the
chain~$\XX$.  Then $\PPbar := \frac{1}{N} \sum_{k=1}^N \PP^k$ is an irreducible
transition matrix with all entries positive.  Moreover, we can effectively use the
original chain to sample from this ``averaged'' chain.  The resulting more general
algorithm (\refA{A:algogen}) obtains, interruptibly, an exact sample~$T$ from the
tree distribution of~$\PP$. 

To describe \refA{A:algogen}, which works in the general setting described at the
outset of \refS{S:algo}, for $t \geq 2 N$ we define the following events to be used
in the context of the algorithm:
\begin{eqnarray*}
A_t &:=& \cap_{l \in [N]} \{X_l(t - 2 N) = 1\}, \\
B_t &:=& \{X_1(t - 2 N + U_0) = 1\}, \\
C_t &:=& \{\{X_1(t - 2 N + U_0 + U_1), X_2(t - 2 N + U_2), \ldots, X_N(t - 2 N +
           U_N)\} = [N]\}, \\
D_t(T)
    &:=& \{\mbox{the graph with directed edges from $X_l(t - 2 N + U_l)$} \\
    &  & \qquad \mbox{to
           $X_l(t - 2 N + U_l + 1)$, $2 \leq l \leq N$, is the arborescence~$T$}\},
           \\
D_t &:=& \cup_{T \in \Tc} D_t(T), \\
E_t(T)
    &:=& A_t \cap B_t \cap C_t \cap D_t(T), \\
E_t &:=& \cup_{T \in \Tc} E_t(T)
      =  A_t \cap B_t \cap C_t \cap D_t.
\end{eqnarray*}
In the following algorithm, successive calls to \Random are assumed to generate
independent random numbers, each uniformly distributed over $[N]$.

\begal[Terminating interruptible stationary sampling]
\label{A:algogen}
{\rm For $t \geq 2 N$, let $E_t(T)$ and $E_t$ be defined as directly above, and let
$E_0 := \emptyset$.  The algorithm is:
\newpage
\begin{code}
\item \s   $t \sset 0$
\item \s   \Repeat
\item \s\s   $t \sset t + 2 N$
\item \s\s   \For $i \sset 0 \To N$
\item \s\s\s   $U_i \sset \Random$
\item \s   \Until $E_t$ holds
\item \s   $S \sset T$, for the unique $T \in \Tc$ such that $E_t(T)$ holds
\item \s   \Return $S$
\end{code}
}
\enal

By
modifying slightly the proof of \refT{T:algoA}, we obtain the following result.
\begin{theorem}
\label{T:algogen}
\refA{A:algogen} is a terminating algorithm for interruptible exact
tree-sampling.
\nopf
\end{theorem}

\begin{remark}
\label{R:theoretical}
Our interest in providing \refA{A:algogen} is more of a theoretical nature
(to establish the possibility of terminating interruptible exact sampling, given
enough copies of a chain) than of a practical nature (to provide an efficient
algorithm).
Thus we have not fine-tuned \refA{A:algogen} to improve its
performance, and we will not analyze its running time here. 
\end{remark}

\begin{remark}
\label{R:nonindependent}
If we make no assumption regarding the independence of the trajectories, then
interruptible sampling becomes impossible for $N \geq 2$ states, no matter how many
trajectories are available.  Indeed, it is then possible that we are in the extreme
case that all the trajectories are identical, \ie, that there is ``really'' only
one trajectory, in which case \refT{T:neg1} applies.
\end{remark}

\section{Impossibility of interruptible exact sampling (I)}
\label{S:neg1}

\refA{A:algogen} requires~$N$ independent synchronized Markov chain trajectories.
This may seem excessive, especially since for interesting chains~$N$ is often
enormously large.  But our next main result, \refT{T:neg1}, shows that this is best
possible.  Note that to prove \refT{T:neg1}, we need only show that interruptible
exact sampling is impossible using $N - 1$ independent trajectories.  Indeed, if
interruptible exact sampling is possible with~$m$ independent trajectories, then
for any $m' \geq m$ it is possible with~$m'$ independent trajectories, since extra
trajectories can always be ignored.

\begin{theorem}
\label{T:neg1}
There
is no algorithm in the passive setting for obtaining an observation from the
stationary distribution of an irreducible aperiodic Markov chain on~$N$ states
which uses fewer than~$N$ independent trajectories from the chain and which is both
interruptible and exact.
\end{theorem}

\begin{proof}
We first establish an equation [\eqref{startsingle}] that must hold if there exists
an interruptible exact sampling algorithm for $N$-state chains (for given $N \geq
2$) that uses only a \emph{single} trajectory; in that case the discussion of
\refS{S:passive} applies verbatim.  A similar equation, namely \eqref{start},
must hold if interruptible exact sampling is possible using $N - 1$ trajectories. 
But~\eqref{start} will lead to a contradiction via a transition-balancing argument.

So we begin with the case of a single trajectory.  Suppose that functions $\phi_{k,
i}$ satisfying (i)--(iii) and (iv$'$) of \refS{S:passive} exist.  We remind the
reader that~(iii) and~(iv$'$) were required to hold for all initial
distributions~$\rho$; throughout the present proof it will suffice to consider
trajectories starting deterministically at~$1$.  Taking~$\rho$ to be unit
mass~$\delta_1$ at~$1$ and $p_{ij}$ to be identically $1 / N$, we find from~(iii)
that, for some $0 \leq k < \infty$,
\begin{equation}
\label{positive}
\sum_{j \in [N]} \sum_{x_1, \ldots, x_k} \phi_{k,j}(1, x_1, \ldots, x_k) > 0
\end{equation}
Let $k_0$ be the minimum such~$k$, and define $\phi_j(\xx) \equiv \phi_j(x_1,
\ldots, x_{k_0}) := \phi_{k_0, j}(1, x_1, \ldots, x_{k_0})$ and $\Xc_j :=
\{\xx = (x_1, \ldots, x_{k_0}): \phi_j(\xx) > 0\}$ for $j \in [N]$.  Again
taking~$\rho$ to be~$\delta_1$ and $p_{ij}$ to be identically $1 / N$, we
find from~(iv$'$) and \eqref{positive} that $\Xc_i \not= \emptyset$ for $i \in
[N]$.  Using~(iv$'$) again, we find that for any transition matrix~$\PP$ with
positive entries and stationary distribution~$\pi$,
\begin{equation}
\label{basic}
\forall i \in [N]:\ \ 
\sum_{\xx \in \Xc_i} \phi_i(\xx)\,p_{1, x_1} \cdots p_{x_{k_0 - 1}, x_{k_0}}
= \pi_i \times \sum_{j \in [N]} \sum_{\xx \in \Xc_j} \phi_j(\xx)\,p_{1, x_1}
\cdots p_{x_{k_0 - 1}, x_{k_0}},
\end{equation}
and all terms on both sides of~\eqref{basic} are positive.  Recalling the notation
of \refS{S:tree}, it now follows in particular that
\begin{equation}
\label{startsingle}
w_2 \sum_{\xx \in \Xc_1} \phi_1(\xx) \prod_{i,j} p^{n_{ij}(\xx)}_{ij} = w_1
\sum_{\xx \in \Xc_2} \phi_2(\xx) \prod_{i,j} p^{n_{ij}(\xx)}_{ij},
\end{equation}
where we write $n_{ij}(\xx)$ for the number of $i \to j$ transitions in the
trajectory $(1, x_1, \ldots, x_{k_0})$ and again all terms on both sides of the
equation are positive.

By the same reasoning, if there exists an interruptible exact sampling algorithm
for $N$-state chains that uses $N - 1$ independent trajectories, then there exist
integer $k \geq 0$ and nonempty sets~$\Xc_1$ and~$\Xc_2$ of $(N - 1)$-tuples
$$
\xx = (x_1(1), \ldots, x_1(k); x_2(1), \ldots, x_2(k); \ldots; x_{N - 1}(1), \ldots,
x_{N - 1}(k))
$$
of $k$-tuples from $[N]$ such that, for any transition matrix~$\PP$ with positive
entries,
\begin{equation}
\label{start}
w_2 \sum_{\xx \in \Xc_1} \phi_1(\xx) \prod_{i,j} p^{n_{ij}(\xx)}_{ij} = w_1
\sum_{\xx \in \Xc_2} \phi_2(\xx) \prod_{i,j} p^{n_{ij}(\xx)}_{ij},
\end{equation}
where, for $l = 1, 2$ and and every $\xx \in \Xc_l$, we have $\phi_l(\xx) > 0$, and
where $n_{ij}(\xx)$ is the sum over $1 \leq m \leq N - 1$ of the numbers of $i \to
j$ transitions within the trajectories $(1, x_m(1), \ldots, x_m(k))$.  To complete
the proof, we will show that~\eqref{start} cannot possibly hold.  We will make key
use of the observation that, for any $\xx \in \Xc_1 \cup \Xc_2$,
\begin{equation}
\label{balance}
0 \leq n_{1+}(\xx) - n_{+1}(\xx) \leq N - 1,
\end{equation} 
where we have introduced the notation
\begin{equation}
\label{ndefs}
n_{1+}(\xx) := \sum_{j = 2}^N n_{1j}(\xx), \qquad n_{+1}(\xx) := \sum_{i = 2}^N
n_{i1}(\xx)
\end{equation}
for the total numbers of transitions out of and into state~$1$, respectively. 
Indeed, since each trajectory $(1, x_m(1), \ldots, x_m(k))$ starts in state~$1$,
the number of transitions out of state~$1$ within such a trajectory either equals
or exceeds by one the number of transitions into state~$1$.

To obtain the desired contradiction, we begin by observing that~\eqref{start} can
be written in the form (eliminating the diagonal variables $p_{ii}$) that
\begin{equation}
\label{starte}
w_2 f_1 = w_1 f_2
\end{equation}
for all $(p_{ij} > 0: 1 \leq i \not= j \leq N)$ such that $\sum_{j: j \neq i}
p_{ij} < 1$ for every $i \in [N]$, where
\begin{equation}
\label{fldef}
f_l := \sum_{\xx \in \Xc_l} \phi_l(\xx) \left[ \prod_{i,j: i \not= j}
p^{n_{ij}(\xx)}_{ij} \right] \left[ \prod_i \Bigl( 1 - \sum_{j: j \not=i} p_{ij}
\Bigr)^{n_{ii}(\xx)} \right], \qquad l = 1, 2.
\end{equation}
Using continuity it follows that~\eqref{starte} holds for all $(p_{ij} \geq 0: 1
\leq i \not= j \leq N)$ such that $\sum_{j: j \neq i} p_{ij} \leq 1$ for every $i
\in [N]$.

For $l = 1, 2$, note that $f_l$ and $w_l$ are both polynomial expressions in the
variables $p_{ij}$, $1 \leq i \not= j \leq N$ (we will denote this entire
collection of $N (N - 1)$ variables by~$\pp$); in fact,
$w_1$ is a polynomial expression in the $(N - 1)^2$ variables $p_{ij}$ with $i, j
\in [N]$ and $i \not\in \{1, j\}$ (with a similar reduction in number of variables
possible for~$w_2$).  Applying \refP{P:identity} (see the Appendix) to $F := w_2
f_1 - w_1 f_2$, we conclude that~\eqref{starte} holds as an equality in the ring of
polynomials in the variables $\pp$ over the complex field.  Henceforth we shall
write $G_1 \equiv G_2$ to indicate such an identity of polynomials~$G_1, G_2$.

According to \refL{L:irred} in the Appendix, the polynomial~$w_1$ (again, over the
complex field) is irreducible; likewise, so is~$w_2$.  From the polynomial
identity $w_2 f_1 \equiv w_1 f_2$ at~\eqref{starte} it then follows that we can
write
\begin{equation}
\label{div}
f_l \equiv w_l f, \qquad l = 1, 2,
\end{equation}
for some polynomial~$f$ in~$\pp$.  Of course, the polynomial identities~\eqref{div}
remain true as we now reduce the number of variables to three by setting $p_{1j}$
to $\alpha$ for $j \not= 1$, $p_{i1}$ to $\beta$ for $i \not= 1$, and $p_{ij}$ to
$\gamma$ if $i \not= 1$, $j \not= 1$, and $i \not= j$.  Observe that now
\begin{eqnarray}
f_l(\alpha, \beta, \gamma)
 &\equiv& \sum_{\xx \in \Xc_l} \phi_l(\xx) \alpha^{n_{1+}(\xx)} \beta^{n_{+1}(\xx)}
            \gamma^{n_{++}(\xx)} \nonumber \\
 &      & \qquad \times [1 - (N - 1) \alpha]^{n_{11}(\xx)} [1 - \beta - (N - 2)
            \gamma]^{H(\xx)},\ \ l = 1, 2,
\label{flabg}
\end{eqnarray}
recalling~\eqref{ndefs} and defining
\begin{equation}
%\label{morendefs}
n_{++}(\xx) := \sum_{i, j \in \{2, \ldots, N\}:\,i \not= j} n_{ij}(\xx), \qquad
H(\xx) := \sum_{i = 2}^N n_{ii}(\xx).
\end{equation}
Also now, by a simple generalization of the bijection argument (\cite{Knuth},
Section~2.3.4.4, p.~390) showing that the number of arborescences rooted at~$1$ is
$N^{N - 2}$,
\begin{equation}
\label{tree-exercise}
w_1(\beta, \gamma) \equiv \beta [\beta + (N - 1) \gamma]^{N - 2};
\end{equation}
and
\begin{equation}
\label{w2}
w_2(\alpha, \beta, \gamma) \equiv \alpha\,v_2(\beta, \gamma)
\end{equation}
for some polynomial $v_2(\beta, \gamma)$ which is not divisible by
$\beta$ [the explanation for~\eqref{w2} being that any $T \in \Tc_2$ has precisely
one directed edge leaving vertex~$1$ and that there exists $T \in \Tc_2$ for
which~$1$ is a leaf.  In fact, it can be shown
%\marginal{sketch:\ we're now in the setting of a reversible chain with stationary
%distribution proportional to $(\beta, \alpha, \alpha, \ldots, \alpha)$}
that $v_2(\beta, \gamma) \equiv [\beta + (N - 1) \gamma]^{N - 2}$, but we won't
need this.]

The idea for the remainder of the proof is to derive from the
identities~\eqref{div}--\eqref{flabg} a polynomial identity in the single
variable~$\beta$, namely~\eqref{gamma-gone}, and then show that~\eqref{gamma-gone}
leads to a contradiction.  We will produce~\eqref{gamma-gone} by eliminating (using
suitable divisibility arguments) first~$\alpha$ and then~$\gamma$.  These arguments
are carried out in the next two lemmas.
\end{proof}

\begin{lemma}
\label{L:gamma-gone}
Suppose that there exists an interruptible exact algorithm in the passive setting
for sampling from the stationary distribution of an irreducible aperiodic Markov
chain on~$N$ states which uses fewer than~$N$ independent trajectories from the
chain.  Then there exist nonempty sets~$\Xc'_1$ and~$\Xc''_1$ and a polynomial~$r$
such that
\begin{equation}
\label{gamma-gone}
\sum_{\xx \in \Xc''_1} \phi_1(\xx) \beta^{n_{+1}(\xx) - m_1(\beta)} (1 -
\beta)^{H(\xx)} \equiv \beta^{N - 2} r(\beta),
\end{equation}
where
$$
m_1(\beta) = \min_{\xx \in \Xc'_1} n_{+1}(\xx).
$$
\end{lemma}

\begin{proof}
Let $m_l(\alpha)$ denote the highest power of~$\alpha$ that divides~$f_l$
at~\eqref{flabg} and define $\mt_l(\alpha) := \min_{\xx \in \Xc_l} n_{1+}(\xx)$.  We
claim that $m_l(\alpha) = \mt_l(\alpha)$, and note that this sort of
highest-power observation will be used frequently---and without accompanying
proof---in the sequel.  [Indeed, $m_l(\alpha) \geq \mt_l(\alpha)$ is clear.  To see
the reverse inequality, divide $f_l$ by $\alpha^{\mt_l(\alpha)}$ and set~$\alpha$
to~$0$ to obtain the expression
\begin{equation}
\label{gl}
\sum_{\xx \in \Xc_l:\,n_{1+}(\xx) = \mt_l(\alpha)} \phi_l(\xx) \beta^{n_{+1}(\xx)}
\gamma^{n_{++}(\xx)} [1 - \beta - (N - 2) \gamma]^{H(\xx)} =: g_l(\beta, \gamma),
\end{equation}
which is not the zero polynomial since it has a positive value when $\beta = 1 / N
= \gamma$.]

By~\eqref{div}, \eqref{tree-exercise}, and~\eqref{w2},
\begin{eqnarray}
m_2(\alpha) = m_1(\alpha) + 1\mbox{\qquad and} \label{malpha} \\
g_l \equiv v_l g, \qquad l = 1, 2, \label{alpha-gone}
\end{eqnarray}
where $g_l$ is the polynomial defined at~\eqref{gl} [recalling $\mt_l(\alpha) =
m_l(\alpha)$], $v_1 := w_1$, $v_2$~is defined at~\eqref{w2}, and $g$ is obtained
from~$f$ by dividing by $\alpha^{m_1(\alpha)}$ and then setting $\alpha = 0$.

Define $\Xc'_l := \{\xx \in \Xc_l: n_{1+}(\xx) = m_l(\alpha)\} \not= \emptyset$
for $l = 1, 2$.  Then, similarly, the highest power $m_l(\beta)$ of~$\beta$
dividing~$g_l$ is $\min_{\xx \in \Xc'_l} n_{+1}(\xx)$;
\begin{equation}
\label{mbeta}
m_2(\beta) = m_1(\beta) - 1;
\end{equation}
and, with
$$
h_1(\beta, \gamma) := \sum_{\xx \in \Xc'_1} \phi_1(\xx) \beta^{n_{+1}(\xx) -
m_1(\beta)} \gamma^{n_{++}(\xx)} [1 - \beta - (N - 2) \gamma]^{H(\xx)},
$$
we have
\begin{equation}
\label{heq}
h_1(\beta, \gamma) \equiv [\beta + (N - 1) \gamma]^{N - 2} h(\beta, \gamma)
\end{equation}
for some polynomial~$h$.

The highest power $m_1(\gamma)$ of~$\gamma$ dividing~$h_1$ is $\min_{\xx \in
\Xc'_1} n_{++}(\xx)$.  Divide both sides of~\eqref{heq} by $\gamma^{m_1(\gamma)}$
and set $\gamma$ to~$0$ to find that~\eqref{gamma-gone} holds for some
polynomial~$r$, where $\Xc''_1 := \{\xx \in \Xc'_1: n_{++}(\xx) = m_1(\gamma)\}
\not= \emptyset$.
\end{proof}

\begin{lemma}
\label{L:gamma-gone-gone}
The identity~\eqref{gamma-gone} cannot hold.
\end{lemma}

\begin{proof}
It follows from~\eqref{gamma-gone} that $n_{+1}(\xx) \geq
m_1(\beta) + N - 2$ for all $\xx \in \Xc''_1$.  But then, for any such~$\xx$ and
some $\xx' \in \Xc'_2$,
\begin{eqnarray*}
n_{+1}(\xx)
 &\geq& m_1(\beta) + N - 2 \\
 &  = & m_2(\beta) + N - 1\mbox{\quad by~\eqref{mbeta}} \\
 &  = & n_{+1}(\xx') + N - 1 \\
 &\geq& n_{1+}(\xx')\mbox{\quad by the second inequality in~\eqref{balance}} \\
 &  = & m_2(\alpha) \\
 &  = & m_1(\alpha) + 1\mbox{\quad by~\eqref{malpha}} \\
 &  = & n_{1+}(\xx) + 1,
\end{eqnarray*}
contradicting the first inequality in~\eqref{balance}.
\end{proof}
\smallskip

\section{Impossibility of interruptible exact sampling (II)}
\label{S:neg2}

\refA{A:algogen} succeeds in using~$N$ independent synchronized Markov chain
trajectories to carry out interruptible exact sampling.  But the algorithm assumes
that each of the trajectories has not only (i)~the same stationary distribution,
but also (ii)~the same transition matrix.  In this section we show (\refT{T:neg2})
that interruptible exact sampling becomes impossible when assumption~(ii) is
dropped, no matter how (finitely) many trajectories are available.

\begin{theorem}
\label{T:neg2}
There is no interruptible algorithm in the passive setting for obtaining an
observation exactly from the common stationary distribution of any finite number of
independent irreducible aperiodic Markov chains on~$N$ states.
\end{theorem}

\begin{proof}
Let~$M$ denote the number of trajectories available.  We first prove the
impossibility of interruptible exact sampling when $M = N = 2$, then more generally
when $N = 2$ (regardless of~$M$), and finally for general~$N$.

For $M = N = 2$, we note that if
\begin{equation}
\label{2cond}
\mbox{$0 < p_{12}, p_{21} < 1$ \quad and \quad $0 < \rho < 1 / \max\{p_{12},
p_{21}\}$},
\end{equation}
then
$$
\PP := \left( \begin{array}{cc}
                1 - p_{12}      & p_{12}     \\
                p_{21}          & 1 - p_{21} \\
              \end{array}
       \right)
\mbox{\qquad and \qquad}
\QQ := \left( \begin{array}{cc}
                1 - \rho p_{12} & \rho p_{12}     \\
                \rho p_{21}     & 1 - \rho p_{21} \\
              \end{array}
       \right)
$$
are irreducible aperiodic transition matrices with common stationary distribution
$$
\pi = \left( \frac{p_{21}}{p_{21} + p_{12}},\ \frac{p_{12}}{p_{21} + p_{12}}
\right).
$$ 

Arguing as in the proof of \refT{T:neg1},
if there exists an interruptible exact sampling algorithm in the present setting,
then there exist integer $k \geq 0$, nonempty sets~$\Zc_1$ and~$\Zc_2$ of pairs
$$
\zz = (\xx, \yy) = (x_1, \ldots, x_k; y_1, \ldots, y_k)
$$
of $k$-tuples from $\{1, 2\}$, and positive numbers $\psi_l(\zz)$ ($\zz \in \Zc_l$,
$l = 1, 2$) such that, whenever~\eqref{2cond} holds,
\begin{equation}
\label{2starte}
p_{12} f_1 = p_{21} f_2
\end{equation}
where, using transition-count notation $n_{ij}$ like that in the proof of
\refT{T:neg1},
\begin{eqnarray}
f_l &=& \sum_{\zz \in \Zc_l} \psi_l(\zz) p^{n_{12}(\zz)}_{12} p^{n_{21}(\zz)}_{21}
          \rho^{n_{12}(\yy) + n_{21}(\yy)} (1 - p_{12})^{n_{11}(\xx)} (1 -
           p_{21})^{n_{22}(\xx)} \nonumber \\
    & & \qquad \quad \times (1 - \rho p_{12})^{n_{11}(\yy)} (1 - \rho
           p_{21})^{n_{22}(\yy)}, \qquad l = 1, 2.
\label{2fldef}
\end{eqnarray}
Using induction on the $\rho$-degree of the polynomial $p_{12} f_1 - p_{21} f_2$ and
\refP{P:identity}, it is easy to show that~\eqref{2starte} holds as an equality in
the ring of polynomials in the variables $p_{12}, p_{21}, \rho$ over the complex
field.
 
For $l = 1, 2$, let
$$
m_l = \min \{[n_{12}(\yy) + n_{21}(\yy)]: \zz = (\xx, \yy) \in \Zc_l\mbox{\ for
some $\xx$}\}
$$ 
denote the highest power of~$\rho$ that divides~$f_l$.  Then,
by~\eqref{2starte}, $m_2 = m_1$.  Divide both sides of~\eqref{2starte} by
$\rho^{m_1}$ and then set $\rho$ to~$0$ to obtain
\begin{equation}
\label{2geqn}
p_{12} g_1 \equiv p_{21} g_2,
\end{equation}
where
$$
g_l := \sum_{\zz \in \Zc'_l} \psi_l(\zz) p^{n_{12}(\zz)}_{12} p^{n_{21}(\zz)}_{21}
(1 - p_{12})^{n_{11}(\xx)} (1 - p_{21})^{n_{22}(\xx)}, \qquad l = 1, 2,
$$
with $\Zc'_l := \{\zz = (\xx, \yy) \in \Zc_l: n_{12}(\yy) + n_{21}(\yy) = m_1\}
\not= \emptyset$.  But [\cf~\eqref{balance} with $N = 2$], if $\zz = (\xx, \yy) \in
\Zc_1 \cup \Zc_2$, then $n_{12}(\yy) = \lc m_1 / 2 \rc$ and $n_{21}(\yy) = \lf m_1 /
2 \rf$.  Dividing both sides of~\eqref{2geqn} by $p^{\lc m_1 / 2 \rc}_{12} p^{\lf
m_1 / 2 \rf}_{21}$ we obtain the polynomial identity
\begin{equation}
\label{2heqn}
p_{12} h_1 \equiv p_{21} h_2,
\end{equation}
where
\begin{equation}
\label{2hldef}
h_l := \sum_{\xx \in \Xc_l} \phi_l(\xx) p^{n_{12}(\xx)}_{12} p^{n_{21}(\xx)}_{21}
(1 - p_{12})^{n_{11}(\xx)} (1 - p_{21})^{n_{22}(\xx)}, \qquad l = 1, 2,
\end{equation}
with
$$
\Xc_l := \{\xx: \mbox{there exists~$\yy$ such that $(\xx, \yy) \in \Zc'_l$}
\} \not= \emptyset, \qquad l = 1, 2
$$
and, for $\xx \in \Xc_l$,
$$
\phi_l(\xx) := \sum_{\yy:\,(\xx, \yy) \in \Zc'_l} \psi_l(\xx, \yy) > 0.
$$

But~\eqref{2heqn} is the case $N = 2$ of~\eqref{start}, which, as shown in the proof
of~\refT{T:neg1}, cannot hold.  This contradiction establishes the theorem in the
case $M = N = 2$.

We leave to the reader the routine extension of the above proof to the case of
arbitrary~$M$ and $N = 2$.  A sketch is that now there are $M - 1$
parameters~$\rho_j$, but by using the same sort of argument for each~$\rho_j$ in
succession that we used above for~$\rho$, one again obtains a contradiction of the
form~\eqref{2heqn} [with $\phi_l(\xx) > 0$ for all $\xx \in \Xc_l \not=
\emptyset$, $l = 1, 2$].

We complete the proof of the theorem by showing that an algorithm for interruptible
exact sampling using~$M$ independent trajectories from chains with $N \geq 3$ states
could be converted into one for two-state chains.

Indeed, while watching independent trajectories of~$M$ generic irreducible
aperiodic two-state chains $X_1, \ldots$, $X_M$ with common (unknown) stationary
distribution $\pi = (\pi_1, \pi_2)$, contemporaneously construct~$M$ independent
irreducible aperiodic $N$-state chains $Y_1, \ldots, Y_M$ by letting
$Y_i(t) = 1$ whenever $X_i(t) = 1$ and selecting an independent uniform random value
from $\{2, \ldots, N\}$ as the value of $Y_i(t)$ at each time~$t$ such that $X_i(t)
= 2$.  The stationary distribution for each~$Y_i$ is $\left(\pi_1, \pi_2 / (N - 1),
\ldots, \pi_2 / (N - 1) \right)$.  Applying the size-$N$ algorithm to $Y_1, \ldots,
Y_M$, suppose the output state is~$S'$.  To finish the construction of the two-state
algorithm, output $S := \min\{S', 2\}$.
\end{proof}

\begin{remark}
\label{R:rev}
(a)~Any two-state chain is reversible, as are the chains~$Y_i$ constructed in the
preceding paragraph.  Thus \refT{T:neg2} remains true even if we assume that the
chains are all reversible.

(b)~Similarly, as mentioned in \refS{S:intro}, interruptible exact sampling from the
stationary distribution is not possible when one observes only a single trajectory
from an irreducible aperiodic reversible finite-state chain. 

(c)~For $N \geq 3$ we do not know whether \refT{T:neg1} remains true if one assumes
that the chain is reversible.
\end{remark}

{\bf Acknowledgment.\ }We thank Dan Naiman for helpful discussions related to the
Appendix.

\def\thesection{A}
\section{Appendix:\ Polynomials}
\label{S:appendix}

In this Appendix we establish two basic facts
about polynomials; these were used in the proof of \refT{T:neg1}.
Throughout the Appendix, we write $F \equiv G$ to indicate that~$F$ and~$G$ are
the same element in the ring of polynomials (in some specified finite collection of
variables) over the complex field.

The first fact is quite simple.  For completeness, we include an elementary
proof. 

\begin{proposition}
\label{P:identity}
Let
$$
\xx = (x_{ij}: 1 \leq i \leq n,\ 1 \leq j \leq k_i)
$$
be a double array of variables, where $n \geq 0$ and $k_i \geq 1$ for $1 \leq i
\leq n$.  If $F(\xx)$ is a polynomial expression that vanishes whenever $x_{ij}
\geq 0$ for all $i, j$ and $\sum_{j = 1}^{k_i} x_{ij} \leq 1$ for all $i$, then
$F(\xx) \equiv 0$. 
\end{proposition}

\begin{proof}
Let $K := \sum_{i = 1}^n k_i$.  The proof is by (strong) induction on $\gk := K +
\deg F$, for which (if~$F$ is not the zero polynomial) the smallest possible value
is $\sum_{1 \leq i \leq 0} 1 + 0 = 0$.  The base case $\gk = 0$ of the induction is
trivial.

For the induction step we may assume $n \geq 1$ and $k_n \geq 1$.  Dividing the
polynomial $F(\xx)$ by $x_{n, k_n}$, we can write
\begin{equation}
\label{division}
F(\xx) \equiv x_{n, k_n} F_1(\xx) + F_2(\xx')
\end{equation}
for polynomials~$F_1$ and~$F_2$, where the variables collection $\xx'$ excludes the
single variable $x_{n, k_n}$.  Setting $x_{n, k_n}$ to~$0$ in~\eqref{division}, we
see that $F_2(\xx')$ is a polynomial satisfying the hypothesis of the proposition;
and (in obvious notation) $K_2 = K - 1$ and $\deg F_2 \leq \deg F$, so that $\gk_2 <
\gk$.  By induction, $F_2(\xx') \equiv 0$, and so from~\eqref{division} we now have
$F(\xx) \equiv x_{n, k_n} F_1(\xx)$.  But now $K_1 = K$ and $\deg F_1 = \deg F -
1$, so that $\gk_1 = \gk - 1$, and one sees that $F_1(\xx)$ satisfies the hypothesis
of the proposition.  By induction, $F_1(\xx) \equiv 0$; we conclude that $F(\xx)
\equiv 0$, as desired.
\end{proof}

As is well known (\eg, \cite{Barshay}, Chapter~4), for any $n \geq 1$ the ring
$\CC[x_1, \ldots, x_n]$ of polynomials in the variables $x_1, \ldots, x_n$ over
the complex field~$\CC$ is a unique factorization domain.  This means that every
nonzero polynomial in $\CC[x_1, \ldots, x_n]$ can be written uniquely (up to complex
scalar multiples) as a (possibly empty) finite product of irreducible
polynomials.  (A polynomial is said to be irreducible if it cannot be factored as
the product of two nonconstant polynomials.)

\begin{lemma}
\label{L:irred}
The polynomial $w_1$ {\rm [}\ie, the polynomial in the $(N - 1)^2$ variables
$p_{ij}$ with $i, j \in [N]$ and $i \notin \{1, j\}$ defined in \refS{S:tree}{\rm
]} is irreducible over the complex field.  
\end{lemma}

\begin{proof}
The proof is by induction on~$N$.  For $N = 1$, the polynomial $w_1 \equiv 1$ (in
no variables) is certainly irreducible.  For $N = 2$, the polynomial $w_1 \equiv
p_{21}$ in the single variable $p_{21}$ is irreducible.  To carry out the induction
step for $N \geq 3$, we will use another induction, on~$l$, to prove the following
claim.
\medskip

\noindent
{\sc Claim.}\ For $3 \leq l \leq N + 1$, let~$y_l$ denote the polynomial in $(N -
1)^2 - (N + 1 - l)$ variables 
obtained from~$w_1$ by setting $p_{m1}$ to~$0$ for $l \leq m
\leq N$.  Then $y_l$ is irreducible for $4 \leq l \leq N + 1$.
\medskip

To prove the claim, we begin by noting that~$y_3$ has the factorization
\begin{equation}
\label{y3fact}
y_3 \equiv p_{21} \omega_2,
\end{equation}
where the polynomial
$$
\omega_2 = \omega_2((p_{ij}: \mbox{$2 \leq i,j \leq N$ and $i \notin \{2, j\}$}))
$$
is obtained from the polynomial~$w_1$ for the state space $[N - 1]$ by changing each
variable name from $p_{ij}$ to $p_{i + 1, j + 1}$.  By the induction hypothesis for
our $N$-induction, $\omega_2$ is irreducible.  Since $p_{21}$ is clearly
irreducible, we conclude that~\eqref{y3fact} is a prime factorization of~$y_3$.

We now treat the base case $l = 4$ of our $l$-induction.  Observe that $y_4
\not\equiv 0$ (consider, e.~g.,\ the tree $N \to N - 1 \to \cdots 2 \to 1$) and
that~$y_4$ is linear in $p_{31}$.  If~$y_4$ is reducible, then we can write
\begin{equation}
\label{y4fact}
y_4 \equiv (g_1 p_{31} + g_2) g_3,
\end{equation}
where~$g_i$ is a polynomial free of the variable $p_{31}$ ($i = 1, 2, 3$) and~$g_3$
is nonconstant.  If we now set $p_{31}$ to~$0$ in~\eqref{y4fact}, the
result is $y_3 \equiv g_2 g_3$.  From the prime factorization~\eqref{y3fact} we
conclude that either $p_{21}$ or~$\omega_2$ divides~$g_3$.  But this is wrong:\
(i)~$p_{21}$ does not divide~$g_3$ because it clearly does not divide~$y_4$
(consider, e.~g.,\ the tree $N \to N - 1 \to \cdots 4 \to 2 \to 3 \to 1$), and
(ii)~$\omega_2$ does not divide~$g_3$ because (we claim) it, too, fails to
divide~$y_4$.  (Indeed, setting
$p_{m2}$ to~$0$ for $3 \leq m \leq N$ causes~$\omega_2$---but clearly
not~$y_4$---to vanish.)  From this contradiction we conclude that~$y_4$ is
irreducible, establishing the $l$-induction base case.

For the $l$-induction step, let $l \geq 5$.  If~$y_l$ is reducible, then we can
write
\begin{equation}
\label{ylfact}
y_l \equiv (h_1 p_{l - 1,1} + h_2) h_3,
\end{equation}
where~$h_i$ is a polynomial free of the variable $p_{l - 1, 1}$ ($i = 1, 2, 3$)
and~$h_3$ is nonconstant.  If we now set $p_{l - 1, 1}$ to~$0$ in~\eqref{ylfact},
the result is $y_{l - 1} \equiv h_2 h_3$.  By the $l$-induction hypothesis, it must
be that~$h_3$ is a nonzero complex scalar multiple of $y_{l - 1}$;
from~\eqref{ylfact} we then deduce that $y_{l - 1}$ divides $y_l$.  But this is
wrong, because setting $p_{m1}$ to~$0$ for $2 \leq m \leq l - 2$ causes $y_{l -
1}$---but clearly not~$y_l$---to vanish.  From this contradiction we conclude
that~$y_l$ is irreducible, completing the $l$-induction.

Finally, set~$l$ to $N + 1$ in the claim to find that $w_1 \equiv y_{N + 1}$ is
irreducible, completing the $N$-induction and the proof of the lemma.
\end{proof}

\end{document}